\documentclass[demo,a4paper]{article}

\usepackage{amsmath}
\usepackage{amssymb}
\usepackage{stmaryrd}
\usepackage{sectsty}
\usepackage{bbm}
\usepackage{blindtext,graphicx}
\usepackage[absolute]{textpos}
\sectionfont{\LARGE}

\title{{\bf On the Riesz dual of ${\bf L^1(\mu)}$}}
\date{}
\author{A. van Rooij}
\begin{document}
	\begin{textblock}{5}(14,1)
	\noindent \Large {\it To Wim Luxemburg}
	\end{textblock}

	\pagenumbering{arabic}
	\maketitle

	\noindent{\bf Abstract} \\
	In this article, $(X,\, \mathcal{A},\, \mu)$ is a measure apace. A classical result establishes a Riesz isomorphism between $L^1(\mu)^{\sim}$ and $L^{\infty}(\mu)$ in case the measure $\mu$ is $\sigma$-finite. In general, there still is a natural Riesz homomorphism $\Phi: L^{\infty}(\mu) \to L^1(\mu)^{\sim},$ but it may not be injective or surjective. We prove that always the range of $\Phi$ is an order dense Riesz subspace of $L^1(\mu)^{\sim}$. If $\mu$ is semi-finite, then $L^1(\mu)^{\sim}$ is a Dedekind completion of $L^{\infty}(\mu)$.

	\paragraph{{\bf Background}}
	Every $h$ in $L^{\infty}(\mu)$ induces an element $\phi_h$ of $L^1(\mu)^{\sim}$ by
		$$\phi_h(f) = \int(fh)d\mu \hspace{1cm} (f \in L^1(\mu)),$$
	As is well known, the map $\Phi: h \mapsto \phi_h$ of $L^{\infty}(\mu)$ into $L^1(\mu)^{\sim}$ is a Riesz isomorphism in case the measure is $\sigma$-finite. Without $\sigma$-finiteness it is still a Riesz homomorphism.

	By Theorem 243 G(a) in [1], $\phi$ is injective if and only if $\mu$ is semi-finite. Given semi-finiteness, $\Phi$ is surjective if and only if $L^{\infty}(\mu)$ is Dedekind complete, as follows from Theorems 241 G(b) and 243 G(b) in [1]. Our result implies that, for semi-finite $\mu ,\, L^1(\mu)^{\sim}$ actually is a Dedekind completion of $L^{\infty}(\mu)$.

	\paragraph{Notations}
	$(X,\, \mathcal{A},\, \mu)$ is a measure space.
	
	We write $L^1,\, L^{\infty},\, \int_X f$ as abbreviations of $L^1(\mu),\, L^{\infty}(\mu),\, \int f d\mu$, respectively.

	A set $A$ in $\mathcal{A}$ determines a measure space $(A,\, \mathcal{A}_A,\, \mu_A)$ where $\mathcal{A}_A = \{Y\cap A: Y\in \mathcal{A}\}=\{Y\in \mathcal{A}: Y \subset A\}$ and $\mu_A$ is the restriction of $\mathcal{A}$ to $\mathcal{A}_A$. We write $L^1(A),\, L^{\infty}(A),\, \int_A f$ for $L^1(\mu_A),\, L^{\infty}(\mu_A),\, \int f d\mu_A,$ respectively.

	A set $A$ in $\mathcal{A}$ is "of $\sigma$-finite measure" is there exist $A_1,\, A_2, \, ...$ in $\mathcal{A}$ with $A = \cup A_n$ and $\mu(A_n) < \infty$ for all $n$.

	If $f \in L^1$, then the (a.e. defined) set $A=\{x: f(x) \neq 0\}$ is of $\sigma$-finite measure and $f=f\mathbbm{1}_A$ (in $L^1$).
	\newpage

	\paragraph{Lemma} {\it The map $h \mapsto \phi_h$ if $L^{\infty}$ into $(L^1)^{\sim}$ is a Riesz homomorphism.}\\
	{\bf Proof} Let $h, \, j \in L^{\infty}$, and let $\alpha \in (L^1)^{\sim}$ be such that $\alpha \geq \phi_h$ and $\alpha \geq \phi_j$ in $(L^1)^{\sim}$; we prove $\alpha \geq \phi_{h \vee j}.$

	If $A:=\{x \in X: h(x) \geq j(x)\}$, then $h \vee j = h \mathbbm{1}_A + j \mathbbm{1}_{X\backslash A}$. For $f \in (L^1)^+$ we see that
		\begin{align*}
		\alpha(f) &= \alpha(f \mathbbm{1}_A) + \alpha(f \mathbbm{1}_{X\backslash A})\\
		&\geq \phi_h(f\mathbbm{1}_A) + \phi_j(f\mathbbm{1}_{X\backslash A})\\
		&= \int_X f\mathbbm{1}_A h + \int_X f\mathbbm{1}_{X\backslash A}j\\
		&= \int_X f(\mathbbm{1}_A h + \mathbbm{1}_{X\backslash A} j) = \phi_{h \vee j}(f).
		\end{align*}

	\paragraph{{\bf Theorem}}{\it The Riesz space $\{ \phi_h : h \in L^{\infty}\}$ is order dense in} $(L^1)^{\sim}$.\\
	{\bf Proof} Let $\alpha \in (L^1)^{\sim}, \, \alpha > 0$; we look for an element $h$ of $L^{\infty}$ with
		$$\hspace{2,94cm}0<\phi_h\leq \alpha.\hspace{2,5cm}(*)$$
	Choose a $g$ in $(L^1)^{+}$ with $\alpha(g) > 0$ and a set $A$ in $\mathcal{A}$ of $\sigma$-finite measure with $g=0$ a.e. off $A$, i.e., $g=g\mathbbm{1}_A$.

	For $u : A \to \mathbb{R}$ define $u^X:X \to \mathbb{R}$ by
		\begin{align*}
		u^X &= u \, \, \, \text{ on } A,\\
		u^X &= 0 \, \, \, \text{ on } X\backslash A.
		\end{align*}
	
	If $u \in L^1(A)$, then $u^X \in L^1$, so $\alpha(u^X)$ exists. By the $\sigma$-finiteness of $A$ there is a  $w \in L^{\infty}(A)$ with $\alpha(u^X) = \int_A uw$ for all $u \in L^1(A)$. Now $h:=w^X$ lies in $L^{\infty}$. For $f$ in $L^1$ we have
		$$\phi_h(f) = \alpha(f\mathbbm{1}_A)$$
since $\phi_h(f)=\int_Xfh = \int_X fw^X = \int_A(f|_A)w=\alpha((f|_A)^X)=\alpha(f\mathbbm{1}_A).$

	Now we can prove the validity of $(*)$:\\
(1) $\phi_h > 0$ because $\phi_h(g) = \alpha(g\mathbbm{1}_A)=\alpha(g)>0$;\\
(2) $\phi_h \leq \alpha$ because for $f \in (L^1)^+$ we have $\phi_h(f) = \alpha(f\mathbbm{1}_A)\leq \alpha(f)$.

	\paragraph{{\bf Corollary}}
	{\it If $\mu$ is semi-finite, then $(L^1)^{\sim}$ is a Dedekind completion of $L^{\infty}$}. (See "Background".)\\

	\noindent As a by-product we obtain the result quoted above: $\Phi$ is surjective if and only if $L^{\infty}$ is Dedekind complete.

	\paragraph{{\bf Reference}} Fremlin, D.H., {\it Measure Theory}, Vol. 2, Torres Fremlin, Colchester (2003).
	


	\paragraph{Adres:} Radboud Universiteit Nijmegen, the Netherlands, Institute for Mathematics, Astrophysics and Particle Physics
	\vspace{0.5 cm}
	
	\noindent maths@math.ru.nl
\end{document}